\documentclass[final]{ws-procs9x6}

\usepackage[latin1]{inputenc}

\usepackage{amsmath}
\usepackage{amssymb}
\usepackage[mathscr]{eucal}

    \usepackage{theorem}

\newenvironment{env}[2]{\begin{#1}#2\end{#1}}{}
    \newcommand{\beq}[1]{\begin{env}{equation}{#1}}
    \newcommand{\beqn}[1]{\begin{env}{equation*}{#1}}
    \newcommand{\bal}[1]{\begin{env}{align}{#1}}
    \newcommand{\baln}[1]{\begin{env}{align*}{#1}}
    \newcommand{\bga}[1]{\begin{env}{gather}{#1}}
    \newcommand{\bgan}[1]{\begin{env}{gather*}{#1}}
    \newcommand{\bflal}[1]{\begin{env}{flalign}{#1}}
    \newcommand{\bflaln}[1]{\begin{env}{flalign*}{#1}}
    \newcommand{\bmu}[1]{\begin{env}{multline}{#1}}
    \newcommand{\bmun}[1]{\begin{env}{multline*}{#1}}
    \newcommand{\bsp}[1]{\begin{env}{split}{#1}}

    \newcommand{\eeq}{\end{env}}
    \newcommand{\eeqn}{\end{env}}
    \newcommand{\eal}{\end{env}}
    \newcommand{\ealn}{\end{env}}
    \newcommand{\ega}{\end{env}}
    \newcommand{\egan}{\end{env}}
    \newcommand{\eflal}{\end{env}}
    \newcommand{\eflaln}{\end{env}}
    \newcommand{\emu}{\end{env}}
    \newcommand{\emun}{\end{env}}
    \newcommand{\esp}{\end{env}}

\renewcommand{\bf}[1]{\textbf{#1}}
\renewcommand{\it}[1]{\textit{#1}}
\renewcommand{\sc}[1]{\textsc{#1}}
\renewcommand{\sf}[1]{\textsf{#1}}
\renewcommand{\sl}[1]{\textsl{#1}}
\renewcommand{\tt}[1]{\texttt{#1}}
\newcommand{\hl}[1]{\bf{\it{#1}}}

\newcommand{\mbf}[1]{\mathbf{#1}}

\newcommand{\cmc}[1]{\mathcal{#1}}
\newcommand{\eus}[1]{\mathscr{#1}}
\newcommand{\euf}[1]{\mathfrak{#1}}
\newcommand{\bb}[1]{\mathbb{#1}}

\newcommand{\nbd}[1]{$#1$\nobreakdash--}
\newcommand{\ol}[1]{\overline{#1}}

\newcommand{\vt}{\vartheta}

\newcommand{\vp}{\varphi}
\newcommand{\om}{\omega}
\newcommand{\Om}{\Omega}

\newcommand{\family}[1]{\left(#1\right)}

\newcommand{\bfam}[1]{\bigl(#1\bigr)}
\newcommand{\Bfam}[1]{\Bigl(#1\Bigr)}
\newcommand{\AB}[1]{\langle#1\rangle}
\newcommand{\bAB}[1]{\bigl\langle#1\bigr\rangle}
\newcommand{\BAB}[1]{\Bigl\langle#1\Bigr\rangle}
\newcommand{\CB}[1]{\{#1\}}
\newcommand{\bCB}[1]{\bigl\{#1\bigr\}}
\newcommand{\BCB}[1]{\Bigl\{#1\Bigr\}}
\newcommand{\SB}[1]{[#1]}

\newcommand{\Matrix}[1]{\begin{pmatrix}#1\end{pmatrix}}

\newcommand{\set}[2][]{
    \ifthenelse{\equal{#1}{}}{
        \CB{#2}}{
        \CB{#1~|~#2}}}
\newcommand{\bset}[2][]{
    \ifthenelse{\equal{#1}{}}{
        \bCB{#2}}{
        \bCB{#1~|~#2}}}
\newcommand{\Bset}[2][]{
    \ifthenelse{\equal{#1}{}}{
        \BCB{#2}}{
        \BCB{#1~\big|~#2}}}

\DeclareMathOperator{\ls}{\normalfont\sf{span}}
\DeclareMathOperator{\cls}{\ol{\ls}}

\DeclareMathOperator{\id}{\normalfont\sf{id}}

\newcommand{\C}{\bb{C}}

\newcommand{\E}{\bb{E}}

\newcommand{\N}{\bb{N}}

\newcommand{\R}{\bb{R}}

\newcommand{\cB}{\cmc{B}}

\newcommand{\sB}{\eus{B}}

\newcommand{\sE}{\eus{E}}
\newcommand{\sF}{\eus{F}}

\newcommand{\sK}{\eus{K}}

\newcommand{\et}{\euf{t}}

\newcommand{\eH}{\euf{H}}

\newcommand{\U}{\mbf{1}}

    \numberwithin{equation}{section}
    \renewcommand{\theequation}{\thesection.\arabic{equation}}

        \newcommand{\definame}{Definition.}
        \newcommand{\propname}{Proposition.}
        \newcommand{\lemname}{Lemma.}
        \newcommand{\exname}{Example.}
        \newcommand{\remname}{Remark.}
        \newcommand{\obname}{Observation.}
        \newcommand{\thmname}{Theorem.}
        \newcommand{\corname}{Corollary.}
        \newcommand{\proofname}{Proof.}
    \theoremheaderfont{\normalfont\bfseries}
    \theoremstyle{change}
        \theorembodyfont{\rmfamily}
            \newtheorem{emp}{}[section]
                \newcommand{\bemp}[1][]{
                    \begin{emp}\hskip-\labelsep\bf{#1}\hskip\labelsep}
                \newcommand{\eemp}{\end{emp}}
\newtheorem{itemp}[emp]{}
                \newcommand{\bitemp}[1][]{
                    \begin{itemp}\hskip-\labelsep\bf{#1}\hskip\labelsep\normalfont\itshape}
                \newcommand{\eitemp}{\end{itemp}}
            \newtheorem{ex}[emp]{\exname}
                \newcommand{\bex}{\begin{ex}}
                \newcommand{\eex}{\end{ex}}
            \newtheorem{defi}[emp]{\definame}
                \newcommand{\bdefi}{\begin{defi}}
                \newcommand{\edefi}{\end{defi}}
            \newtheorem{rem}[emp]{\remname}
                \newcommand{\brem}{\begin{rem}}
                \newcommand{\erem}{\end{rem}}
            \newtheorem{ob}[emp]{\obname}
                \newcommand{\bob}{\begin{ob}}
                \newcommand{\eob}{\end{ob}}
        \theorembodyfont{\normalfont\itshape}
            \newtheorem{thm}[emp]{\thmname}
                \newcommand{\bthm}{\begin{thm}}
                \newcommand{\ethm}{\end{thm}}
            \newtheorem{Prop}[emp]{\propname}
                \newcommand{\bprop}{\begin{Prop}}
                \newcommand{\eprop}{\end{Prop}}
            \newtheorem{cor}[emp]{\corname}
                \newcommand{\bcor}{\begin{cor}}
                \newcommand{\ecor}{\end{cor}}
            \newtheorem{lem}[emp]{\lemname}
                \newcommand{\blem}{\begin{lem}}
                \newcommand{\elem}{\end{lem}}
    \newcommand{\qedsymbol}{~\rule[-0.35mm]{2mm}{2mm}}
    \newcounter{PROOF}[emp]
    \newenvironment{Proof}[1]{
        \vspace{1ex}
        \renewcommand{\item}[1][\stepcounter{PROOF}(\roman{PROOF})]%
            {##1\hskip\labelsep}
        \noindent\textsc{#1\hskip\labelsep}}{
        \nolinebreak\qedsymbol}
    \newcommand{\PROOF}[1][\proofname]{
        \begin{Proof}{#1}\ignorespaces}
    \newcommand{\qed}{\end{Proof}}
    \newcommand{\noqed}{
        \renewcommand{\qedsymbol}{}
        \end{Proof}}

\makeatletter
\def\cite{\@cghitrue\@ifnextchar [{\@tempswatrue
    \@citex}{\@tempswatrue\@citex[]}}
\def\@cite#1#2{{[{#1}]\if@tempswa \fi}}
\def\@citex[#1]#2{\if@filesw\immediate\write\@auxout{\string\citation{#2}}\fi
 \@tempcnta\z@\@tempcntb\m@ne\def\@citea{}\@cite{\@for\@citeb:=#2\do
  {\@ifundefined
   {b@\@citeb}{\@citeo\@tempcntb\m@ne\@citea\def\@citea{}%
     {\mbox{\bfseries ?}}\@warning
   {Citation `\@citeb' on page \thepage \space undefined}}%
  {\setbox\z@\hbox{\global\@tempcntc0\csname b@\@citeb\endcsname\relax}%
   \ifnum\@tempcntc=\z@ \@citeo\@tempcntb\m@ne
    \@citea\def\@citea{,}\hbox{\csname b@\@citeb\endcsname}%
   \else
    \advance\@tempcntb\@ne
    \ifnum\@tempcntb=\@tempcntc
    \else\advance\@tempcntb\m@ne\@citeo
    \@tempcnta\@tempcntc\@tempcntb\@tempcntc\fi\fi}}\@citeo}{#1}}
\makeatother

\begin{document}

\renewcommand{\tt}[1]{\small\texttt{#1}}
%\renewcommand{\hl}[1]{{\footnotesize\slshape\sf{#1}}}

%\addtolength{\hoffset}{6mm}
%\addtolength{\textheight}{5mm}
%\setlength{\evensidemargin}{\oddsidemargin}

%\bibliographystyle{alpha}

\title{A Problem of Powers and the Product
\\
of Spatial Product Systems\footnote{This work is supported by a PPP-project by DAAD and DST and by a Research in Pairs project at MfO.}}

\author{B.V.\ Rajarama Bhat}
\address{Statistics and Mathematics Unit, Indian Statistical Institute Bangalore, R.\ V.\ College Post, Bangalore 560059, India, E-mail: \tt{bhat@isibang.ac.in}, Homepage: \tt{http://www.isibang.ac.in/Smubang/BHAT/}}

\author{Volkmar Liebscher}
\address{Institut für Mathematik und Informatik, Ernst-Moritz-Arndt-Universität Greifswald, 17487 Greifswald, Germany, E-mail: \tt{volkmar.liebscher@uni-greifswald.de}, Homepage: \tt{http://www.math-inf.uni-greifswald.de/biomathematik/liebscher/}}

\author{Michael Skeide\footnote{MS is supported by research funds of the University of Molise and the Italian MIUR (PRIN 2005) and by research funds of the Dipartimento S.E.G.e S.\ of University of Molise.}}
\address{Dipartimento S.E.G.e S., Università degli Studi del Molise, Via de Sanctis, 86100 Campobasso, Italy, E-mail: \tt{skeide@math.tu-cottbus.de}, Homepage: \tt{http://www.math.tu-cottbus.de/INSTITUT/lswas/\_skeide.html}}

%\renewcommand{\baselinestretch}{1}

%\vspace{2ex}

\vfill

\begin{abstract}
In the 2002 AMS summer conference on ``Advances in Quantum Dynamics'' in Mount Holyoke Robert Powers proposed a sum operation for spatial \nbd{E_0}semigroups. Still during the conference Skeide showed that the Arveson system of that sum is the product of spatial Arveson systems. This product may but need not coincide with the tensor product of Arveson systems. The Powers sum of two spatial \nbd{E_0}semigroups is, therefore, up to cocycle conjugacy Skeide's product of spatial noises.
\end{abstract}

%\tableofcontents
%}

%\clearpage

%\vfill

\keywords{Quantum dynamics, quantum probability, Hilbert modules, product systems, $E_0$--semigroups. 2000 AMS-Subject classification: 46L53; 46L55; 46L08; 60J25; 12H20.}

%\tableofcontents

%\vfill

%{\parskip0.5ex plus 0.5ex minus 0.5ex

\bodymatter
\renewcommand{\thefootnote}{[\arabic{footnote}]}
    \numberwithin{equation}{section}
    \renewcommand{\theequation}{\thesection.\arabic{equation}}

\newpage

\section{Introduction} \label{intro}

Let $\vt^i=\bfam{\vt^i_t}_{t\in\R_+}$ $(i=1,2)$ be two \nbd{E_0}semigroups on $\sB(H)$ with associated Arveson systems ${\eH^i}^\otimes=\bfam{\eH^i_t}_{t\in\R_+}$ (Arveson \cite{Arv89}). Furthermore, let $\Om^i=\bfam{\Om^i_t}_{t\in\R_+}\subset\sB(H)$ be two semigroups of intertwining isometries for $\vt^i$ (units). Then
\beqn{\label{*}\tag{$*$}
T_t
\Matrix{
a_{11}&a_{12}\\a_{21}&a_{22}
}
~=~
\Matrix{
\vt^1_t(a_{11})&\Om^1_ta_{12}{\Om^2_t}^*\\\Om^2_ta_{12}{\Om^1_t}^*&\vt^2_t(a_{22})
}
}\eeqn
defines a CP-semigroup $T=\bfam{T_t}_{t\in\R_+}$ on $\sB(H\oplus H)$. In the 2002 AMS summer conference on ``Advances in Quantum Dynamics'' Robert Powers asked for the Arveson system associated with $T$ (Bhat \cite{Bha96}, Arveson \cite{Arv97}). During that conference (see \cite{Ske03c}) Skeide showed that this product system is nothing but the \it{product of the spatial product systems} introduced in Skeide \cite{Ske06d} (published first in 2001). Meanwhile, Powers has formalized the above \it{sum operation} in \cite{Pow04} and he has proved that the product may, but need not coincide with the tensor product of the involved Arveson systems, a fact suspected already in \cite{Ske03c,Lie00p1}.

In these notes we extend Powers' construction to the case of spatial \nbd{E_0}semigroups $\vt^i$ on $\sB^a(E^i)$ where $E^i$ are Hilbert \nbd{\cB}modules. We obtain the same result as in \cite{Ske03c}, namely, the product system of the minimal dilation of the CP-semigroup on $\sB^a(E^1\oplus E^2)$ defined in analogy with \eqref{*} is the product of spatial product systems from \cite{Ske06d}.

Like in \cite{Ske03c}, it is crucial to understand the following point (on which we will spend some time in Section \ref{CPB(E)}): In Bhat and Skeide \cite{BhSk00} to every CP-semigroup on a \nbd{C^*}algebra $\cB$ a product system of \nbd{\cB}algebra has been constructed. However, the \nbd{C^*}algebra in question here is $\sB^a(E^1\oplus E^2)$, where $\sB^a(E)$ denotes the algebra of all adjointable operators on a Hilbert module $E$. So what has the product system of \nbd{\sB^a(E^1\oplus E^2)}correspondences to do with the product systems of the \nbd{E_0}semigroups $\vt^1$ and $\vt^2$, which are product systems of \nbd{\cB}correspondences? The answer to this question, like in \cite{Ske03c}, will allow to construct the Arveson system of a CP-semigroup on $\sB(H)$ \it{without} having to find first its minimal dilation. To understand even the Hilbert space case already requires, however, module techniques.

\section{Product systems, CP-semigroups, $E_0$--semigroups and dilations}

Throughout these notes, by $\cB$  we denote a unital \nbd{C^*}algebras. There are no spatial product systems where $\cB$ is nonunital. There is no reasonable notion of unit for product systems of correspondences over nonunital \nbd{C^*}algebras, where $\cB$ could not easily be substituted by a unital ideal of $\cB$.

\bemp[Product systems.~]
Product systems of Hilbert modules (\hl{product system} for short) occurred in different contexts in Bhat and Skeide \cite{BhSk00}, Skeide \cite{Ske02}, Muhly and Solel \cite{MuSo02} and other more recent publications. Let $\cB$ be a unital \nbd{C^*}algebra. A product system is a family $E^\odot=\bfam{E_t}_{t\in\R_+}$ of \hl{correspondences} $E_t$ over $\cB$ (that is, a (right) Hilbert \nbd{\cB}module with a unital representation of $\cB$) with an associative identification
\beqn{
E_s\odot E_t
~=~
E_{s+t},
}\eeqn
where $E_0=\cB$ and for $s=0$ or $t=0$ we get the canonical identifications. By $\odot$ we denote the (internal) tensor product of correspondences.

If we want to emphasize that we do not put any technical condition, we say \hl{algebraic} product system. There are concise definitions of \hl{continuous} \cite{Ske03b} and \hl{measurable} (separable!) \cite{Hir04} product systems of \nbd{C^*}correspondences, and \hl{measurable} (separable pre-dual!) \cite{MuSo07} product systems of \nbd{W^*}correspondences. Skeide \cite{Ske07p} will discuss \hl{strongly continuous} product systems of von Neumann correspondences. We do not consider such constraints in these notes. We just mention for the worried reader that the result from \cite{Ske06d,Ske03b} that the product of continuous spatial product systems is continuous.
\eemp

\bemp[Units.~]
A \hl{unit} for a product system $E^\odot$ is family $\xi^\odot=\bfam{\xi_t}_{t\in\R_+}$ of elements $\xi_t\in E_t$ such that
\beqn{
\xi_s\odot\xi_t
~=~
\xi_{s+t}
}\eeqn
and $\xi_0=\U\in\cB=E_0$. A unit may be \hl{unital} ($\AB{\xi_t,\xi_t}=\U\forall t\in\R_+$), \hl{contractive} ($\AB{\xi_t,\xi_t}\le\U\forall t\in\R_+$), or \hl{central} ($b\xi_t=\xi_tb\forall t\in\R_+,b\in\cB$).

We do not pose technical conditions on the unit. But, sufficiently continuous units can be used to pose technical conditions on the product system in a nice way; see \cite{Ske03c}.
\eemp

\bemp[The product system of a CP-semigroup.~]\label{CPPS}
Let $T=\bfam{T_t}_{t\in\R_+}$ be a (not necessarily unital) CP-semigroup on a unital \nbd{C^*}algebra $\cB$. According to Bhat and Skeide \cite{BhSk00} there exists a product system $E^\odot$ with a unit $\xi^\odot$ determined uniquely up to isomorphism (of the pair $(E^\odot,\xi^\odot)$) by the following properties:
\begin{enumerate}
\item
$\AB{\xi_t,b\xi_t}=T_t(b)$.

\item
$E^\odot$ is \hl{generated} by $\xi^\odot$, that is, the smallest subsystem of $E^\odot$ containing $\xi^\odot$ is $E^\odot$.
\end{enumerate}
In analogy with Paschke's \cite{Pas73} GNS-construction for CP-maps, we call $(E^\odot,\xi^\odot)$ the \hl{GNS-system} of $T$ and we call $\xi^\odot$ the \hl{cyclic unit}. In fact, $\sE_t=\cls\cB\xi_t\cB$ is the \hl{GNS-module} of $T_t$ with cyclic vector $\xi_t$. For the comparison of the product system of Powers' CP-semigroup with a product of product systems it is important to note that
\bal{\label{GNSind}
E_t
&
~=~
\cls\bCB{x^n_{t_n}\odot\ldots\odot x^1_{t_1}\colon n\in\N,t_n+\ldots+t_1=t,x^k_{t_k}\in\sE_{t_k}}
\\\notag
&
~=~
\cls\bCB{b_n\xi_{t_n}\odot\ldots\odot b_1\xi_{t_1}b_0\colon n\in\N,t_n+\ldots+t_1=t,b_k\in\cB}.
}\eal
In fact, the product system $E_t$ can be obtained as an inductive limit of the expressions $\sE_{t_n}\odot\ldots\odot \sE_{t_1}$ over refinement of the partitions $t_n+\ldots+t_1=t$ of $\SB{0,t}$.
\eemp

\bemp[The product system of an $E_0$--semigroup on $\sB^a(E)$.~]\label{E0PS}
Let $E$ be a Hilbert \nbd{\cB}module with a \hl{unit vector} $\xi$ (that is, $\AB{\xi,\xi}=\U$) and let $\vt=\bfam{\vt_t}_{t\in\R_+}$ be an \hl{\nbd{E_0}semigroup} (that is, a semigroup of unital endomorphisms) on $\sB^a(E)$. Let us denote by $xy^*$ $(x,y\in E)$ the \hl{rank-one operator}
\beqn{
xy^*
\colon
z
~\longmapsto~
x\AB{y,z}.
}\eeqn
Then $p_t:=\vt_t(\xi\xi^*)$ is a projection and the range $E_t:=p_tE$ is a Hilbert \nbd{\cB}submodule of $E$. By defining the (unital!) left action $bx_t=\vt_t(\xi b\xi^*)x_t$ we turn $E_t$ into a \nbd{\cB}correspondence. One easily checks that
\beqn{
x\odot y_t
~\longmapsto~
\vt_t(x\xi^*)y_t
}\eeqn
defines an isometry $u_t\colon E\odot E_t\rightarrow E$. Clearly, if $\vt_t$ is \hl{strict} (that is, precisely, if $\cls\vt_t(EE^*)E=E$), then $u_t$ is a unitary. Identifying $E=E\odot E_t$ and using the semigroup property, we find
\bal{\label{assvtdef}
\vt_t(a)
~&=~
a\odot\id_{E_t}
&
(E\odot E_s)\odot E_t
~&=~
E\odot(E_s\odot E_t).
}\eal
The restriction of $u_t$ to $E_s\odot E_t$ is a bilinear unitary onto $E_{s+t}$ and the preceding associativity reads now $(E_r\odot E_s)\odot E_t=E_r\odot(E_s\odot E_t)$. Obviously, $E_0=\cB$ and the identifications $E_t\odot E_0=E_t=E_0\odot E_t$ are the canonical ones. Thus, $E^\odot=\bfam{E_t}_{t\in\R_+}$ is a product system.

For \nbd{E_0}semigroups on $\sB(H)$ the preceding construction is due to Bhat \cite{Bha96}, the extension to Hilbert modules to Skeide \cite{Ske02}. We would like to mention that Bhat's construction does not give the Arveson system of an \nbd{E_0}semigroup, but its \it{opposite} Arveson system (all orders in tensor products reversed). By Tsirelson \cite{Tsi00p1} the two need not be isomorphic. For Hilbert \nbd{C^*}modules Arveson's construction does not work. For von Neumann modules it works, but gives a product system of von Neumann correspondences over $\cB'$, the commutant of $\cB$; see \cite{Ske03c,Ske04p'}.

Existence of a unit vector is not a too hard requirement, as long as $\cB$ is unital. (If $E$ has no unit vector, then a finite multiple $E^n$ will have one; see \cite{Ske04p'}. And product systems do not change under taking direct sums.) We would like to mention a further method to construct the product system of an \nbd{E_0}semigroup, that works also for nonunital $\cB$. It relies on the representations theory of $\sB^a(E)$ in Muhly, Skeide and Solel \cite{MSS06}. See \cite{Ske04p'} for details.
\eemp

\bemp[Dilation and minimal dilation.~]\label{dmdSEC}
Suppose $E^\odot$ is a product system with a unit $\xi^\odot$. Clearly, $T_t:=\AB{\xi_t,\bullet\xi_t}$ defines a CP-semigroup $T=\bfam{T_t}_{t\in\R_+}$, which is unital, if and only if $\xi^\odot$ is unital. Obviously, $E^\odot$ is the product systems of $T$, if and only if it is generated by $\xi^\odot$.

If $\xi^\odot$ is unital, then we may embed $E_t$ as $\xi_s\odot E_t$ into $E_{s+t}$. This gives rise to an inductive limit $E$ and a factorization $E=E\odot E_t$, fulfilling the associativity condition in \eqref{assvtdef}. It follows that $\vt_t(a)=a\odot\id_{E_t}$ defines an \nbd{E_0}semigroup $\vt=\bfam{\vt_t}_{t\in\R_+}$ on $\sB^a(E)$. The embedding $E_t\rightarrow E_{s+t}$ is, in general, only right linear so that, in general, $E$ is only a right Hilbert module.

Under the inductive limit all $\xi_t\in E_t\subset E$ correspond to the same unit vector $\xi\in E$. Moreover, $\xi=\xi\odot\xi_t$, so that the vector expectation $\vp:=\AB{\xi,\bullet\xi}$ fulfills $\vp\circ\vt_t(\xi b\xi^*)=T_t(b)$, that is, $(E,\vt,\xi)$ is a \hl{weak dilation} of $T$ in the sense of \cite{BhPa94,BhSk00}. Clearly, the product system of $\vt$ (constructed with the unit vector $\xi$) is $E^\odot$.

Suppose $\vt$ is a strict \nbd{E_0}semigroup on some $\sB^a(E)$ and that $\xi$ is a unit vector in $E$. One may show (see \cite{Ske02}) that $T_t(b):=\vp\circ\vt_t(\xi b\xi^*)$ defines a (necessarily unital) CP-semigroup (which it dilates), if and only if the projections $p_t:=\vt_t(\xi\xi^*)$ increase. In this case, the product system $E^\odot$ of $\vt$ has a unit $\xi^\odot=\bfam{\xi_t}_{t\in\R_+}$ with $\xi_t:=p_t\xi$, which fulfills $T_t=\AB{\xi_t,\bullet\xi_t}$. We say the weak dilation $(E,\vt,\xi)$ of $T$ is \hl{minimal}, if the \hl{flow} $j_t(b):=\vt_t(\xi b\xi^*)$ generates $E$ out of $\xi$. One may show that this is the case, if and only if the product system of $\vt$ coincides with the product system of $T$. The minimal (weak) dilation is determined up to suitable unitary equivalence.
\eemp

\brem
We would like to emphasize that in order to construct the minimal dilation of a unital CP-semigroup $T$, we first constructed the product system of $T$ and then constructed the dilating \nbd{E_0}semigroup $\vt$ (giving back the product system of $T$). It is not necessary to pass through minimal dilation to obtain the product system of $T$, but rather the other way round.
\erem

\bemp[Spatial product systems.~]
Following \cite{Ske06d}, we call a product system $E^\odot$ \hl{spatial}, if it has central unital \hl{reference unit} $\om^\odot=\bfam{\om_t}_{t\in\R_+}$. The choice of the reference unit is part of the spatial structure, so we will write a pair $(E^\odot,\om^\odot)$. For instance, a \hl{morphism} $w^\odot\colon E^\odot\rightarrow F^\odot$ between product systems $E^\odot$ and $F^\odot$ is a family $w^\odot=\bfam{w_t}_{t\in\R_+}$ of mappings $w_t\in\sB^{a,bil}(E_t,F_t)$ (that is, bilinear adjointable mappings from $E_t$ to $F_t$) fulfilling $w_s\odot w_t=w_{s+t}$ and $w_0=\id_\cB$. To be a \hl{spatial} morphism of spatial product systems, $w^\odot$ must send the reference unit of $E^\odot$ to the reference unit of $F^\odot$.

Our definition matches that of Powers \cite{Pow87} in that an \nbd{E_0}semigroup $\vt$ on $\sB^a(E)$ admits a so-called intertwining semigroup of isometries, if and only if the product system of $\vt$ is spatial. It does not match the usual definition for Arveson systems, where an Arveson system is \it{spatial}, if it has a unit. The principle result of Barreto, Bhat, Liebscher and Skeide \cite{BBLS04} asserts that a product system of von Neumann correspondences is spatial, if it has a (continuous) unit. But, for Hilbert modules this statement fails. In fact, we show in \cite{BLS08p} that, unlike for Arveson systems, a subsystem of a product system of Fock modules need not be spatial (in particular, it need not be Fock).

There are many interesting questions about spatial product systems, open even in the  Hilbert space case. Does the spatial structure of the spatial product system depend on the choice of the reference unit? The equivalent question is, whether every spatial product system is \it{amenable} \cite{Bha01} in the sense that the product system automorphisms act transitively on the set of units. Tsirelson \cite{Tsi04p1} claims they are not. But, still there is a gap that has not yet been filled. In contrast to this, the question raised by Powers [Pow04], whether the product defined in the next section depends on the reference units, or not, we can answer in the negative sense; see \cite{BLS08p1}.
\eemp

\bemp[The product of spatial product systems.~]
The basic motivation of \cite{Ske06d} was to define an \it{index} of a product system and to find a \it{product} of product systems under which the index is \it{additive}. Both problems could not be solved in full generality, but precisely for the category of spatial product systems.

The mentioned result \cite{BLS08p} is one of the reasons why it is hopeless to define an index for nonspatial product systems. However, once accepted the necessity to restrict to spatial product systems (anyway, the index of a nonspatial Arveson systems is somewhat an artificial definition), everything works as we know it from Arveson systems, provided we indicate the good product operation.

In the theory of Arveson systems, there is the tensor product (of arbitrary Arveson systems). However, for modules this does not work. (You may write down the tensor product of correspondences, but, in general, it is not possible to define a product system structure.)

The \hl{product} of two spatial product systems ${E^i}^\odot$ $(i=1,2)$ with reference units ${\om^i}^\odot$ is the spatial product system $(E^1\circledcirc E^2)^\odot$ with reference unit $\om^\odot$ which is characterized uniquely up to spatial isomorphism by the following properties:
\begin{enumerate}
\item\label{subid}
There are spatial isomorphisms ${w^i}^\odot$ from ${E^i}^\odot$ onto subsystems of $(E^1\circledcirc E^2)^\odot$.

\item\label{gen}
$(E^1\circledcirc E^2)^\odot$ is generated by these two subsystems.

\item\label{reforth}
$\AB{w^1_t(x^1_t),w^2_t(y^2_t)}=\AB{x^1_t,\om^1_t}\AB{\om^2_t,y^2_t}$.
\end{enumerate}
Existence of the product follows by an inductive limit; see \cite{Ske06d}. By Condition \ref{subid} we may and, usually, will identify the factors as subsystem of the product. Condition \ref{reforth} means, roughly speaking, that the reference units of the two factors are identified, while components from different factors which are orthogonal to the respective reference unit are orthogonal in the product. Condition \ref{gen} means that
\bmun{
(E^1\circledcirc E^2)_t
\\
~=~
\cls\bCB{x^n_{t_n}\odot\ldots\odot x^1_{t_1}\colon n\in\N,t_n+\ldots+t_1=t,x^k_{t_k}\in E^i_{t_k}~(i=1,2)}.
}\emun
It is important to note (crucial exercise!) that this may be rewritten in the form
\bmu{\label{prodgen}
(E^1\circledcirc E^2)_t
\\
~=~
\cls\bCB{x^n_{t_n}\odot\ldots\odot x^1_{t_1}\colon n\in\N,t_n+\ldots+t_1=t,x^k_{t_k}\in\sE^i_{t_k}~(i=1,2)},
}\emu
where we put $\sE_t:=\cB\om_t\oplus(E^1_t\ominus\cB\om^1_t)\oplus(E^2_t\ominus\cB\om^2_t)$ (the direct sum of $E^1$ and $E^2$ with ``identification of the reference vectors'' and denoting the new reference vector by $\om_t$). Written in that way, it is easy to see that the subspaces are actually increasing of the partitions $t_n+\ldots+t_1=t$ of $\SB{0,t}$. This gives an idea how to obtain the product as an inductive limit; see \cite{Ske06d}.
\eemp

\section{The product system of $\cB$--correspondences of a CP-semigroup on $\sB^a(E)$}\label{CPB(E)}

In Section \ref{CPPS} we have said what the product system of CP-semigroup on $\cB$ is. It is a product systems of \nbd{\cB}correspondences. On the other hand, if \nbd{\sB(H)}people speak about the product system of a unital CP-semigroup on $\sB(H)$, they mean an Arveson system, that is, a product system of Hilbert spaces. Following Bhat \cite{Bha96} and Arveson \cite{Arv97}, the Arveson system of a unital CP-semigroup is the Arveson system of its minimal dilating \nbd{E_0}semigroup. (To be specific, we mean the product system constructed as in Section \ref{E0PS} following \cite{Bha96}, not the product system constructed in \cite{Arv89}, which is anti-isomorphic to the former.) A precise understanding of the relation between the two product systems, one of \nbd{\sB(H)}modules, the other of Hilbert spaces, will allow to avoid the construction of the minimal dilation. But we will discuss it immediately for CP-semigroups on $\sB^a(E)$.

Suppose we have a Hilbert \nbd{\sB^a(E)}module $F$. Then we may define the Hilbert \nbd{\cB}module $F\odot E$. Every $y\in F$ gives rise to a mapping $y\odot\id\in\sB^a(E,F\odot E)$ defined by $(y\odot\id_E)x=y\odot x$ with adjoint $y^*\odot\id_E\colon y'\odot x\mapsto\AB{y,y'}x$. These mappings fulfill $(y\odot\id_E)^*(y'\odot\id_E)=\AB{y,y'}$ and $ya\odot\id_E=(y\odot\id_E)a$ for every $a\in\sB^a(E)$. Via $a\mapsto y\odot\id_E$ we may identify $F$ as a subset of $\sB^a(E,F\odot E)$. This subset is strictly dense but, in general, it need not coincide. In fact, we have always $F\supset\sK(E,F\odot E)$ where the \hl{compact operators} between Hilbert \nbd{\cB}modules $E_1$ and $E_2$ are defined as $\sK(E_1,E_2):=\cls\CB{x_1x_2^*\colon x_i\in E_i}$, and $F=\sK(E,F\odot E)$ whenever the right multiplication is strict (in the same sense as left multiplication, namely, $\cls FEE^*=F$).

\brem\label{stritprem}
The space $\sB^a(E,F\odot E)$ may be thought of as the strict completion of $F$, and  it is possible to define a strict tensor product of \nbd{\sB^a(E)}correspondences. We do not need this here, and refer the interested reader to \cite{Ske04p'}.
\erem

Now suppose that $F$ is a \nbd{\sB^a(E)}correspondence with strict left action. If $E$ has a unit vector, then, doing as in Section \ref{E0PS}, we see that $F$ factors into $E\odot F_E$ (where the $F_E$ is a suitable multiplicity correspondence from $\cB$ to $\sB^a(E)$) and $a\in\sB^a(E)$ acts on $F=E\odot F_E$ as $a\odot\id_{F_E}$. For several reasons we do not follow Section \ref{E0PS}, but refer to the representation theory of $\sB^a(E)$ from \cite{MSS06}. This representation theory tells us that $F_E$ may be chosen as $E^*\odot F$, where $E^*$ is the \hl{dual} \nbd{\cB}\nbd{\sB^a(E)}correspondence of $E$ with operations $\AB{x^*,x'^*}:=xx'^*$ and $bx^*a:=(a^*xb^*)^*$. Then, clearly,
\bmun{
F
~=~
\cls\sK(E)F
~=~
\sK(E)\odot F
\\
~=~
(E\odot E^*)\odot F
~=~
E\odot(E^* \odot F)
~=~
E\odot F_E
}\emun
explains both how the isomorphism is to be defined and what the action of $a$ is. Putting this together with the preceding construction, we obtain
\beqn{
\sB^a(E,E\odot E_F)
~\supset~
F
~\supset~
\sK(E,E\odot E_F)
~=~
E\odot E_F\odot E^*,
}\eeqn
where we defined the \nbd{\cB}correspondence $E_F:=E^*\odot F\odot E$.

\brem
We do not necessarily have equality $F=\sK(E,E\odot E_F)$. But if we have (so that $F$ is a \it{full} Hilbert \nbd{\sK(E)}module), then the operation of \it{tensor conjugation} with $E^*$ may be viewed as an operation of Morita equivalence for correspondences in the sense of Muhly and Solel \cite{MuSo00}. In what follows, the generalization to Morita equivalence of product systems \cite{Ske04p'} is in the background. An elaborate version for the strict tensor product (see Remark \ref{stritprem}) can be found in \cite{Ske04p'}.
\erem

We observe that the assignment (the functor, actually) $F\mapsto E_F:=E^*\odot F\odot E$ respects tensor products. Indeed, if $F_1$ and $F_2$ are \nbd{\sB^a(E)}correspondences with strict left actions, then
\bmu{\label{MePSdef}
E_{F_1}\odot E_{F_2}
~=~
(E^*\odot F_{F_1}\odot E)\odot(E^*\odot F_{F_2}\odot E)
\\
~=~
E^*\odot F_{F_1}\odot(E\odot E^*\odot F_{F_2})\odot E
\\
~=~
E^*\odot F_{F_1}\odot F_{F_2}\odot E
~=~
E_{F_1\odot F_2}.
}\emu
It is, clearly, associative. It respects inclusions and, therefore, inductive limits. If $E$ is full, then $E^*\odot\sB^a(E)\odot E=\cB$. We summarize:

\bitemp[Proposition \cite{Ske04p'}.~]\label{FEprop}
Suppose that $E$ is full (for instance, $E$ has a unit vector). Suppose that $F^\odot=\bfam{F_t}_{t\in\R_+}$ is a product system of \nbd{\sB^a(E)}correspondences such that the left actions of all $F_t$ are strict.

Then the family $E^\odot=\bfam{E_t}_{t\in\R_+}$ of \nbd{\cB}correspondences $E_t:=E^*\odot F_t\odot E$ with product system structure defined by \eqref{MePSdef} is a product system.

Moreover, if the $F_t$ are inductive limits over families $\sF_\et$, then the $E_t$ are inductive limits over the corresponding $\sE_\et:=E^*\odot\sF_\et\odot E$.
\eitemp

\bthm
Let $F^\odot$ be the GNS-system of a strict unital CP-semigroup $T$ on $\sB^a(E)$, and denote by $(F,\theta,\zeta)$ the minimal dilation of $T$.

Then $E^\odot$ (from Proposition \ref{FEprop}) is the product system of  the strict \nbd{E_0}semigroup $\vt$ induced on $\sB^a(F\odot E)\cong\sB^a(F)\odot\id_E=\sB^a(F)$ by $\theta$.

The triple $(F\odot E,\vt,p=\vt_0(\zeta\zeta^*))$ is the unique minimal dilation of $T$ to the operators on a Hilbert \nbd{\cB}module in the sense that
\beqn{
p(F\odot E)
~=~
\U_{\sB^a(E)}\odot E
~=~
E
}\eeqn
and
\beqn{
p\vt_t(a)p
~=~
T_t(a).
}\eeqn
\ethm

\noindent\bf{Proof.~}
We proceed precisely as in the proof of [\refcite{Ske04p'}, Theorem 5.12]. We know (see Section \ref{dmdSEC}) that the product system of the minimal $\theta$ is $F^\odot$. Though, we have a unit vector $\zeta$ in $F$, it is more suggestive to think of the correspondences $F_t$ to be obtained as $F_t=F^*\odot{_t}F$ where $_tF$ is $F$ viewed as \nbd{\sB^a(E)}correspondences via $\theta_t$; see [\refcite{Ske04p'}, Section 2] for details. In the same way, the product system of $\vt$ is $(F\odot E)^*\odot{_t}(F\odot E)$. We find
\bmun{
(F\odot E)^*\odot{_t}(F\odot E)
~=~
(E^*\odot F^*)\odot({_t}F\odot E)
\\
~=~
E^*\odot(F^*\odot{_t}F)\odot E
~=~
E^*\odot F_t\odot E
~=~
E_t.
}\emun
(Note: The first step where $_t$ goes from outside the brackets into, is just the definition of $\vt_t$.) This shows the first statement.

For the second statement, we observe that $x\mapsto\zeta\odot x$ provides an isometric embedding of $E$ into $F\odot E$ and that $p$ is the projection on the range $\zeta\odot E$ of this embedding. Clearly,
\bmun{
p\vt_t(a)p
~=~
(\zeta\zeta^*\odot\id_E)(\theta_t(a)\odot\id_E)(\zeta\zeta^*\odot\id_E)
\\
~=~
\Bfam{\zeta\bAB{\zeta,\theta_t(a)\zeta}\zeta^*}\odot\id_E
~=~
(\zeta T_t(a)\zeta^*)\odot\id_E
~=~
T_t(a),
}\emun
when $\sB^a(E)$ is identified with the corner $(\zeta\sB^a(E)\zeta^*)\odot\id_E$ in $\sB^a(F\odot E)$.\qedsymbol

\brem
If $T$ is a normal unital CP-semigroup on $\sB(H)$ (normal CP-maps on $\sB(H)$ are strict), then $E^\odot$ is nothing but the Arveson system of $T$ (in the sense of Bhat's construction). Note that we did construct  $E^\odot$ \bf{without} constructing the minimal dilation first. In the theorem the minimal dilation occurred only, because we wanted to verify that our product system coincides with the one constructed via minimal dilation.
\erem

\brem
We hope that the whole discussion could help to clarify the discrepancy between the terminology and constructions in the case of CP-semigroups on $\sB(H)$ and those for CP-semigroups on $\cB$. The semigroups of this section lie in between, in that they are CP-semigroups on $\sB^a(E)$, so not general $\cB$ but also not just $\sB(H)$. The operation that transforms the product system of \nbd{\sB^a(E)}correspondences into a product system of \nbd{\cB}correspondences is \it{cum grano salis} an operation of Morita equivalence. (In the von Neumann case and when $E$ is full, it is Morita equivalence.) We obtain \nbd{\cB}correspondences because $E$ is a Hilbert \nbd{\cB}module. For $\sB(H)$ we obtain \nbd{\C}correspondences (or Hilbert spaces), because $H$ is a Hilbert \nbd{\C}module.
\erem

\section{Powers' CP-semigroup}\label{PCP}

We, finally, come to  Powers' CP-semigroup and to the generalization to Hilbert modules of the result from \cite{Ske03c} that its product system is the product of the involves spatial product systems.

Let $\vt^i$ $(i=1,2)$ be two strict \nbd{E_0}semigroup on $\sB^a(E^i)$ ($E^i$ two Hilbert \nbd{\cB}modules with unit vectors $\om^i$) with spatial product systems ${E^i}^\odot$ (as in Section \ref{E0PS}) and unital central reference units ${\om^i}^\odot$. Since $\om^i_t$ commutes with $\cB$, the mapping $b\mapsto\om^i_tb$ is bilinear. Consequently, $\Om^i_t:=\id_{E^i}\odot\om^i_t\colon x^i\mapsto x^i\odot\om^i_t\in E^i\odot E^i_t=E^i$ defines a semigroup of isometries in $\sB^a(E^i)$. (The isometries are \hl{intertwining} in the sense that $\vt^i_t(a)\Om^i_t=(a\odot\id_{E^i_t})(\id_{E^i}\odot\om^i_t)=(\id_{E^i}\odot\om^i_t)a=\Om^i_ta$.) It follows that
\beqn{
T_t\Matrix{a_{11}&a_{12}\\a_{21}&a_{22}}
~=~
\Matrix{\vt^1_t(a_{11})&\Om^1_ta_{12}{\Om^2_t}^*\\\Om^2_ta_{21}{\Om^1_t}^*&\vt^2_t(a_{22})}
}\eeqn
defines a unital semigroup on $\sB^a\bfam{\text{\raisebox{.1ex}{$\substack{E^1\\E^2}$}}}$. (We see later on that $T_t$ is completely positive, by giving its GNS-module explicitly.) Using the identifications $E^i=E^i\odot E^i_t$ we find the more convenient form
\beqn{
T_t\Matrix{a_{11}&a_{12}\\a_{21}&a_{22}}
~=~
\Matrix{a_{11}\odot\id_{E^1_t}&(\id_{E^1}\odot\om^1_t)a_{12}(\id_{E^2}\odot{\om^2_t}^*)\\(\id_{E^2}\odot\om^2_t)a_{21}(\id_{E^1}\odot{\om^1_t}^*)&a_{22}\odot\id_{E^i_t}}
}\eeqn
where $T_t$ maps from $\sB^a\bfam{\text{\raisebox{.1ex}{$\substack{E^1\\E^2}$}}}$ to $\sB^a\bfam{\text{\raisebox{0ex}{$\substack{E^1\odot E^1_t\\E^2\odot E^2_t}$}}}=\sB^a\bfam{\text{\raisebox{.1ex}{$\substack{E^1\\E^2}$}}}$.

Denote by $F^\odot$ the product system of $T$ in the sense of Section \ref{E0PS}, that is, the $F_t$ are \nbd{\sB^a\bfam{\text{\raisebox{.1ex}{$\substack{E^1\\E^2}$}}}}correspondences. By Proposition \ref{FEprop}, setting $E_t:=\bfam{\text{\raisebox{.1ex}{$\substack{E^1\\E^2}$}}}^*\odot F_t\odot\bfam{\text{\raisebox{.1ex}{$\substack{E^1\\E^2}$}}}$ we define a product system $E^\odot$ of \nbd{\cB}correspondences and
\beqn{
\sK\family{\bfam{\text{\raisebox{.1ex}{$\substack{E^1\\E^2}$}}},\bfam{\text{\raisebox{.1ex}{$\substack{E^1\\E^2}$}}}\odot E_t}
~\subset~
F_t
~\subset~
\sB^a\family{\bfam{\text{\raisebox{.1ex}{$\substack{E^1\\E^2}$}}},\bfam{\text{\raisebox{.1ex}{$\substack{E^1\\E^2}$}}}\odot E_t}.
}\eeqn

\bthm
$E^\odot$ is the product $(E^1\circledcirc E^2)^\odot$ of the spatial product systems ${E^1}^\odot$ and ${E^2}^\odot$.
\ethm

\noindent\bf{Proof.~}
Recall that, by \eqref{GNSind}, $F_t$ is the inductive limit of expressions of the form
\beqn{
\sF_\et
~:=~
\sF_{t_n}\odot\ldots\odot\sF_{t_1}
}\eeqn
over the partitions $\et=(t_n,\ldots,t_1)$ with $t_n+\ldots+t_1=t$, where $\sF_t$ is the GNS-module of $T_t$ with cyclic vector $\zeta_t$.

Put $\sE_t=\bfam{\text{\raisebox{.1ex}{$\substack{E^1\\E^2}$}}}^*\odot\sF_t\odot\bfam{\text{\raisebox{.1ex}{$\substack{E^1\\E^2}$}}}$. Then $\sF_t\subset\sB^a\family{\bfam{\text{\raisebox{.1ex}{$\substack{E^1\\E^2}$}}},\bfam{\text{\raisebox{.1ex}{$\substack{E^1\\E^2}$}}}\odot\sE_t}$. We claim that $\sE_t=\cB\om_t\oplus(E^1_t\ominus\cB\om^1_t)\oplus(E^2_t\ominus\cB\om^2_t)$ and that $\zeta_t$ is the operator given by
\beqn{
\zeta_t
\Matrix{
x^1\\x^2
}
~=~
\Matrix{z^1\\0}\odot\bfam{\AB{\om^1_t,y^1_t}\,,\,p^1_ty^1_t\,,\,0}
+
\Matrix{0\\z^2}\odot\bfam{\AB{\om^2_t,y^2_t}\,,\,0\,,\,p^2_ty^2_t},
}\eeqn
with $\bfam{\text{\raisebox{.1ex}{$\substack{x^1\\x^2}$}}}=\bfam{\text{\raisebox{.1ex}{$\substack{z^1\odot y^1_t\\z^2\odot y^2_t}$}}}\in\bfam{\text{\raisebox{.1ex}{$\substack{E^1\\E^2}$}}}=\bfam{\raisebox{.1ex}{\text{$\substack{E^1\odot E^1_t\\E^2\odot E^2_t}$}}}$ and $p^i_t:=\id_{E^i_t}-\om^i_t{\om^i_t}^*$. To show this, we must check two things. Firstly, we must check whether $\AB{\zeta_t,a\zeta_t}=T_t(a)$. This is straightforward and we leave it as an exercise. Secondly, we must check whether elements of the form on the right-hand side of
\beqn{
\Matrix{
x^1\\x^2
}^*
\odot
\zeta_t
\odot
\Matrix{z^1\odot y^1_t\\z^2\odot y^2_t}
~\longmapsto~
\family{
\Matrix{
x^1\\x^2
}^*
\odot
\id_{\sE_t}
}
\zeta_t
\Matrix{z^1\odot y^1_t\\z^2\odot y^2_t}
}\eeqn
are total in $\sE_t$. For the right-hand side we find
\bmun{
\BAB{\Matrix{x^1\\x^2},\Matrix{z_1\\0}}\bfam{\AB{\om^1_t,y^1_t}\,,\,p^1_ty^1_t\,,\,0}
+
\BAB{\Matrix{x^1\\x^2},\Matrix{0\\z_2}}\bfam{\AB{\om^2_t,y^2_t}\,,\,0\,,\,p^2_ty^2_t}
\\
~=~
\bfam{\AB{x_1,z_1}\AB{\om^1_t,y^1_t}+\AB{x_2,z_2}\AB{\om^2_t,y^2_t}\,,\,\AB{x_1,z_1}p^1_ty^1_t\,,\,\AB{x_2,z_2}p^2_ty^2_t}.
}\emun
From this, totality follows.

By Proposition \ref{FEprop}, we obtain $E_t$ as inductive limit over the expressions
\beqn{
\sE_\et
~:=~
\sE_{t_n}\odot\ldots\odot\sE_{t_1},
}\eeqn
which, by the preceding computation, precisely coincides with what is needed, according to \eqref{prodgen}, to obtain $E_t=(E^1\circledcirc E^2)_t$.\qedsymbol

%\bibliography{mybib}
\newcommand{\Swap}[2]{#2#1}\newcommand{\Sort}[1]{}

%\listofOWs

\end{document}